\documentclass[12pt, dvipdfmx]{amsart}
\usepackage[utf8]{inputenc}
\usepackage{textcomp}
\usepackage{mathtools, amsmath, amssymb, mathrsfs}
\usepackage{tikz}
\usepackage{url, hyperref}
\usepackage{amsthm}
\usepackage{comment}
\usepackage{kantlipsum}

\calclayout

\theoremstyle{definition}
\newtheorem{counter}{counter}[section]
\newtheorem{definition}[counter]{Definition}
\newtheorem{theorem}[counter]{Theorem}
\newtheorem{proposition}[counter]{Proposition}
\newtheorem{lemma}[counter]{Lemma}

\newtheorem*{Proof}{Proof}

\newtheorem{remark}[counter]{Remark}
\newtheorem{claim}{Claim}

\makeatletter
\newcommand{\tpitchfork}{
  \vbox{
    \baselineskip\z@skip
    \lineskip-.52ex
    \lineskiplimit\maxdimen
    \m@th
    \ialign{##\crcr\hidewidth\smash{$-$}\hidewidth\crcr$\pitchfork$\crcr}
  }
}
\makeatother

\renewcommand{\qed}{\blacksquare}

\newcommand{\set}[2]{\left\{ #1 \ \middle| \ #2 \right\}}

\newcommand{\xto}{\xrightarrow}

\newcommand{\isom}{\cong}

\newcommand{\toto}{\rightrightarrows}

\renewcommand{\le}{\leqslant}
\renewcommand{\ge}{\geqslant}
\renewcommand{\phi}{\varphi}
\renewcommand{\epsilon}{\varepsilon}

\newcommand{\hy}{{\rm \mathchar`-}}
\renewcommand{\tilde}{\widetilde}
\newcommand{\void}{\varnothing}

\newcommand{\op}{\mathrm}
\renewcommand{\frak}{\mathfrak}
\newcommand{\cal}{\mathcal}
\renewcommand{\bold}{\mathbb}

\newcommand{\R}{\mathbb R}

\newcommand{\X}{\mathfrak X}

\title{Canonical stratification of definable Lie groupoids
}
\author{Masato TANABE}
\address{(Masato Tanabe) M2, Graduate School of Information Science and Technology, Hokkaido University, Sapporo 060-0814, Japan}
\email{tanabe.masato.i8@elms.hokudai.ac.jp}
\date{\today}

\begin{document}

\begin{abstract}
Our aim is to precisely present a tame topology counterpart to canonical stratification of a Lie groupoid.
We consider a  definable Lie groupoid in semialgebraic, subanalytic, o-minimal over $\R$, or more generally, Shiota's $\X$-category. 
We show that there exists a canonical Whitney stratification of the Lie groupoid into definable strata which are invariant under the groupoid action.
This is a generalization and refinement of results on real algebraic group action which J.~N.~Mather and V.~A.~Vassiliev independently stated with sketchy proofs.
A crucial change to their proofs is to use Shiota's isotopy lemma and approximation theorem in the context of tame topology. 
\end{abstract}

\maketitle

\tableofcontents

\setcounter{section}{0}

\section{Introduction}

It is a basic problem to find some nice decomposition of a given space, e.g., triangulation, cellular decomposition, and Whitney stratification.
As a particular example, while less known, J.~N.~Mather \cite[Theorem 1]{Mat_IDGA} and V.~A.~Vassiliev \cite[Theorem 8.6.6]{Vass} independently have shown that for a real algebraic manifold with real algebraic group action, there exists a $C^\omega$ Whitney stratification into invariant semialgebraic strata.
In their works, major examples are jet spaces $J^r(n,p)$ with actions of algebraic groups $\cal{A}^r$, $\cal{K}^r$, $\cal{R}^r$ of $r$-jets of diffeomorphism-germs, which take an important role in singularity theory.
Also in a different context, the existence of invariant Whitney stratification of a proper Lie groupoid in $C^\infty$ category is very recently discussed in Crainic--Mestre \cite{Cra}.
That is also related to a new trend in homotopy theory on conically smooth stratifications (cf. Lurie \cite{Lur} and Ayala--Francis--Tanaka \cite{AFT}).

The result of Mather or Vassiliev has potential to be generalized in two directions.
The one is the extension from group actions to groupoids as in \cite{Cra}, and the other is from the semialgebraic category to an o-minimal category over $\R$ (van den Dries \cite{vdD}) or $\X$-category (Shiota \cite{Shi_G}). 
An o-minimal category or $\X$-category is an axiomatic generalization of semialgebraic category and subanalytic category, respectively, and the notion of definable Lie groupoid in o-minimal category can be found in, e.g., \cite{Hru}. 

Now we state a semialgebraic version of the main theorem. 

\begin{theorem}\label{main}
{\em Let $\cal{G} \toto M$ be a semialgebraic $C^\omega$ (that is Nash) Lie groupoid.
Then, there exists a filtration
\[M = M_0 \supset M_1 \supset M_2 \supset \cdots \supset M_{d+1} = \void\]
of $M$ such that for each $i = 0, 1, \dots, d \, (= \dim M)$,}
\begin{enumerate}
\item {\em the set $M_i$ is a $\cal{G}$-invariant semialgebraic closed subset of $M$;}
\item {\em the set $M_i - M_{i+1}$ is a semialgebraic $C^\omega$ manifold of codimension $i$ in $M$ (unless it is empty) and $\{M_i - M_{i+1}\}_{i=0}^d$ is a Whitney stratification of $M$;}
\item {\em the quotient space $(M_i - M_{i+1}) / \cal{G}$ admits a $C^\omega$ manifold structure and the quotient map $q \colon M_i - M_{i+1} \to (M_i - M_{i+1}) / \cal{G}$ is a $C^\omega$ locally trivial fibration. Moreover, the quotient manifold and the quotient map are piecewise algebraic.}
\end{enumerate}
\end{theorem}

Here are some remarks.
In (2), we can see that $\cal{G}|_{M_i - M_{i+1}} \toto M_i - M_{i+1}$ is a regular Lie groupoid on each connected component of $M_i - M_{i+1}$.
In (3), {\em piecewise algebraic spaces and maps} are the notions introduced by Kontsevich--Soibelman \cite[Appendix]{Kon} (roughly, such a space is made by gluing semialgebraic subsets via semialgebraic isomorphisms). Note that $(M_i - M_{i+1}) / \cal{G}$ may have several connected components of different dimension. 
Applying this theorem to real algebraic manifold with real algebraic group action $\cal{G} = G\times M \toto M$, we recover the result of Mather and Vassiliev (precisely saying, our claim is a bit stronger, because we also show the piecewise algebraicity in (3).

The above theorem will be presented in a slightly different form for a general o-minimal category over $\R$, and also for subanalytic or general $\X$-category with the assumption that both $\cal{G}$ and $M$ are bounded (Theorem \ref{main2}).  
That is, the word `semialgebraic' in the statement above is replaced by `definable' within a more general context. However, we have to restrict the regularity of manifolds and maps to be of class $C^r$ with $1 \le r < \infty$, instead of $C^\omega$. The reason is due to a key tool used in our proof. 

The proofs of Theorem \ref{main} and \ref{main2} are close to Mather's idea in \cite{Mat_IDGA}.
A crucial difference from Mather's is the following point. 
To show (3), Mather used the original version of Thom's isotopy lemma \cite{Mat_NTS}, which is achieved by integrating a stratified vector field ($C^\infty$ on strata) and approximating the obtained $C^\infty$ trivialization map by a $C^\omega$ one, so the finally obtained trivialization map is not semialgebraic in general (also Vassiliev referred to Artin's theory on algebraic spaces for finding $C^\omega$ charts in (3)).
That does not fit with the context of tame topology, namely, we have to avoid to use integration or algebraic spaces. 
Instead of that, we employ, as specially refined techniques, 
Shiota's isotopy lemma and approximation theorem within an o-minimal or $\X$-category \cite[(II.5.2), (II.6.1)]{Shi_G}. This approximation theorem is only shown for the class of $C^r$ with $0 \le r < \infty$ in general, i.e., it is unsolved yet for $r = \infty, \omega$. 
Thus Theorem \ref{main2} is stated only with the regularity of finite order.
On the other hand, the approximation theorem in semialgebraic and subanalytic category holds under any regularity, i.e., $r$ can be $\infty$, $\omega$ \cite{Shi_N, Whi_0}.  
Thus above Theorem \ref{main} is stated with the regularity of class $C^\omega$ in (2) and (3). 
Note again that our proof uses an essentially different tool from that of Mather and Vassiliev.

The present paper consists of the following sections. 
In \S 2, we recall some notions of o-minimal category over $\R$ and $\X$-category, and Whitney stratifications and Lie groupoids in the category. We also recall two key tools, isotopy lemma and approximation theorem in the category.
In \S 3, we state and prove our main theorem.

\section{Preliminaries}

In this section, we recall some notions of geometry of definable sets and maps, and a definable version of stratification theory and Lie groupoid theory. We also recall two key tools, isotopy lemma and approximation theorem for definable maps.
Our main geometric objects in the present paper are {\em definable Lie groupoids}.

Hereafter, $r$ denotes a fixed integer such that $1 \le r < \infty$.

\subsection{Definable sets and maps}

We firstly recall definitions of {\em Shiota's $\X$ and $\X_0$-category} or an {\em o-minimal category} over the real field $\R$.
This is a legitimate generalization of semialgebraic and subanalytic category.
For a detail, see \cite{Shi_G, vdD}.

\begin{definition}[$\X$ and $\X_0$-category {\cite[Chap.II, p.95, p.146]{Shi_G}}]
$ $\\
(1) Let $\X$ be a family of subsets of Euclidean spaces which satisfies the following axioms:
\begin{description}
\item[Axiom (i)] Every algebraic set in any Euclidean space is an element of $\X$;
\item[Axiom (ii)] If $X_1, X_2$ are elements of $\X$, then $X_1 \cap X_2$, $X_1 - X_2$, and $X_1 \times X_2$ are elements of $\X$;
\item[Axiom (iii)] If $X \subset \R^n$ is an element of $\X$ and $p \colon \R^n \to \R^{n-1}$ is a linear map such that the restriction to the closure of $X$ in $\R^n$ is proper, then $p(X) \in \X$;
\item[Axiom (iv)] If $X \subset \R$ is an element of $\X$, then each point of $X$ has a neighborhood in $X$ which is a finite union of points and intervals.
\end{description}
(2) Let $\X_0$ be a family $\X$ which satisfies the following stronger axioms than the above (iii) and (iv):
\begin{description}
\item[Axiom (iii)$_0$] If $X \subset \R^n$ is an element of $\X_0$, and $p \colon \R^n \to \R^{n-1}$ is a linear map, then $p(X) \in \X_0$;
\item[Axiom (iv)$_0$] If $X \subset \R$ is an element of $\X_0$, then $X$ is a finite union of points and intervals.
\end{description}
An {\em $\X$-set} is an element of $\X$, and an {\em $\X$-map} is a continuous map between $\X$-sets whose graph is an $\X$-set. Similary we define an {\em $\X_0$-set} and an {\em $\X_0$-map}.
An $\X_0$-category is the same as an {\em o-minimal category over $\R$}.
\end{definition}

\begin{remark}[boundedness condition for $\X$-maps] \label{boundedness-conditions}
Let $X' \subset X \subset \R^m$ and $Y' \subset Y \subset \R^n$ be $\X_0$-sets and $f \colon X \to Y$ an $\X_0$-map.
Thanks to the axiom (iii)$_0$,  the preimage $f^{-1}(Y')$ and the image $f(X')$ are $\X_0$-sets. However, for an $\X$-map $f \colon X \to Y$  between $\X$-sets, the same claim does not hold in general. Then we often require the following boundedness condition for $f$:
\begin{enumerate}
\item[(I)] For any bounded subset $B \subset \R^m$, $f(X \cap B) \subset \R^p$ is bounded.
\item[(II)] For any bounded subset $C \subset \R^n$, $f^{-1}(C) \subset \R^n$ is bounded.
\end{enumerate}
According to \cite[(II.1.1), (II.1.6)]{Shi_G}, the boundedness condition ensures that $f^{-1}(Y')$ and $f(X')$ are $\X$-sets, and that for another $\X$-map $g \colon Y \to \R^p$, $g \circ f$ is an $\X$-map. 
\end{remark}

Throughout the present paper, we work over such a category (semialgebraic, subanalytic, $\X$-category, and o-minimal category over $\R$), and use the words, {\em definable} sets/maps, for simplicity.

\subsection{Definable Whitney stratifications}
We recall some notions and facts of stratification theory.
In this paper, a definable set $X \subset \R^m$ is said to be a ($k$-dimensional) {\em definable $C^r$ manifold} if $X$ is a ($k$-dimensional) $C^r$ regular submanifold of $\R^m$ (that is, every point $x \in X$ has a neighborhood $U$ in $\R^m$ and a definable $C^r$ diffeomorphism $\phi \colon U \to \phi(U) \subset \R^m$ such that $\phi(X \cap U) = \phi(U) \cap (\R^k \times \{0\})$). 
Let $V \subset \R^m$ denote a definable set.

\begin{definition}[definable $C^r$ stratification]
A {\em definable $C^r$ stratification} of $V$ is a partition $\cal{S}$ of $V$ into definable $C^r$ manifolds which is locally finite at each point of $\R^m$, i.e., 
for each $x \in \R^m$, there exists a neighborhood $U$ of $x$ in $\R^m$ such that $\# \{S \in \cal{S} \ | \ S \cap U \ne \void\} < \infty$.
Each member belonging to $\cal{S}$ is called {\em stratum} of $\cal{S}$.
\end{definition}

\begin{proposition}[{\cite[(II.1.8)]{Shi_G}}]
{\em Every definable set admits a definable $C^r$ stratification.}
\end{proposition}

We define the {\em dimension} of $V$ as the highest dimension of strata belonging to a stratification, and write it $\dim V$.

\begin{definition}[regular/singular point]
Let $x \in V$ be a point.
A point $x \in V$ is said to be {\em $C^r$ regular} of dimension $k$ if there exists an open neighborhood $U$ of $x$ in $\R^m$ such that $V \cap U$ is a $k$-dimensional $C^r$ regular submanifold of $\R^m$.
A point $x \in V$ is said to be {\em $C^r$ singular} if $x$ is not regular of highest dimension (i.e., $\dim V$). 
\end{definition}

Let $\Sigma_r V$ denote the set of all $C^r$ singular points of $V$, called the {\em $C^r$ singular set} of $V$.

\begin{lemma}[{\cite[(II.1.10)]{Shi_G}}]\label{sing-set}
{\em The set $\Sigma_r V$ is a definable closed subset of $V$ of dimension less than $\dim V$.}
\end{lemma}

Let $X_\alpha, X_\beta \subset \R^m$ be definable $C^r$ manifolds.
We here omit the definition of {\em the Whitney (b)-regularity condition of $X_\alpha$ over $X_\beta$ at a point $y \in X_\beta$}.
See \cite{Gib_TSSM, Mat_NTS, Shi_G} for a detail.


\begin{definition}[bad set]
The {\em bad set $B(X_\alpha, X_\beta)$ of the pair $(X_\alpha, X_\beta)$} is the subset of $X_\beta$ consisting of points at which $X_\alpha$ fails to be Whitney (b)-regular over $X_\beta$.
\end{definition}

\begin{lemma}[{\cite[(II.1.13)]{Shi_G}}, {\cite[Lemma 2.4]{Tro_GPEDWS}}]\label{bad-set}
{\em The set $B(X, Y)$ is a definable subset of $Y$ of dimension less than $\dim Y$.}
\end{lemma}

\begin{definition}[the frontier condition]
The pair $(X_\alpha, X_\beta)$ satisfies the {\em frontier condition} if $\op{Cl}_{\R^m} X_\alpha \cap X_\beta \ne \void$ implies that $\op{Cl}_{\R^m} X_\alpha \supset X_\beta$, where $\op{Cl}_A B$ means the closure of $B$ in $A$. 
\end{definition}

\begin{definition}[definable Whitney stratification] 
Let $V \subset \R^m$ be a definable set. 
A definable $C^r$ stratification $\cal{S}$ of $V$ is said to be {\em definable $C^r$ Whitney stratification} if every pair of strata of $\cal{S}$ satisfies the Whitney (b)-regularity condition and the frontier condition.
\end{definition}

We know the existence theorem of definable $C^r$ Whitney stratification of definable subsets  \cite[(2.7)]{Gib_TSSM}, \cite[(I.2.2), (II.1.14)]{Shi_G}. For later use in our proof of the main result, we state it in a slightly modified form.

\begin{proposition}\label{classical-stratification}
{\em Let $V \subset \R^n$ be a definable closed subset and $V' \subset V$ a definable $C^r$ manifold such that $\dim V = \dim V' > \dim (V - V')$.
Then, $V$ admits a finite definable $C^r$ Whitney stratification such that $V'$ is the top stratum.}
\end{proposition}

\begin{Proof}
We put $d = \dim V = \dim V'$ and define a filtration of $V$ by closed subsets inductively:
\[V = V_0 \supset V_1 \supset \cdots \supset V_d \supset V_{d+1} = \void,\]
where $V_0 \coloneqq V$, $V_1 \coloneqq V - V'$, and for each $i \ge 1$,
\[V_{i+1} \coloneqq
\begin{cases}
V_i & (\op{codim} \ V_i > i) \\
\op{Cl}_V \left( \Sigma_r V_i \cup B(V', V_i - \Sigma_r V_i) \cup \bigcup_{j=1}^{i-1} B(V_j - \Sigma_r V_j, V_i - \Sigma_r V_i) \right) & (\op{codim} \ V_i = i).
\end{cases}
\]
By Lemma \ref{sing-set} and Lemma \ref{bad-set}, we can obtain that $\op{codim} \ V_i \ge i$ and each $V_i$ is definable, thus the induction works.
From this construction, the partition $\{V_i - V_{i-1}\}_{i=0}^d$ of $V$ is a definable $C^r$ Whitney stratification. $\qed$
\end{Proof}

\begin{remark}[regularity ($r = \omega$) in semialgebraic/subanalytic case]\label{remark1}
Lemma \ref{sing-set} holds even if we replace ``definable $C^r$'' as ``semialgebraic $C^\omega$ (that is Nash)'' and ``subanalytic $C^\omega$'' \cite[Proposition 9.7.4]{Michel-et-al}, \cite[Theorem 1.2.2 (v)]{Tamm}. Thus, Proposition \ref{classical-stratification} also holds.
\end{remark}

\subsection{Isotopy lemma and Approximation theorem}

We state a definable version of isotopy lemma and approximation theorem proved by Shiota in \cite{Shi_G}. 

\subsubsection{Isotopy lemma}
The definable version of isotopy lemma will critically be used in the proof of our main result, especially, in finding a fibration structure of the quotient map and finding a manifold structure of the quotient space. See \cite{Mat_NTS} for the original version --- Thom's isotopy lemma.

\begin{definition}[definable $C^r$ locally trivial fibration]
Let $X \subset \R^n$ be a definable set and $Y \subset \R^p$ a definable $C^r$ manifold.
\begin{enumerate}
\item A definable map $f \colon X \to Y$ is said to be {\em of class $C^r$} if $f$ is extended to a $C^r$ map from a neighborhood of $X$ in $\R^n$ to $\R^p$.
\item A definable $C^r$ map $f \colon X \to Y$ is said to be a {\em definable $C^r$ locally trivial fibration} if for each point $y \in Y$, there exist a neighborhood $V$, which is definable, of $y$ in $Y$ and a definable $C^r$ map $\Phi \colon f^{-1}(V) \to f^{-1}(y)$ such that the map $(f, \Phi) \colon f^{-1}(V) \to V \times f^{-1}(y)$ is $C^r$ diffeomorphic.
\end{enumerate}
The map $(f, \Phi)$ is called a {\em trivialization map} of $f$ on $V$.
\end{definition}

\begin{theorem}[Shiota's (first) isotopy lemma {\cite[(II.6.1)]{Shi_G}}] \label{isotopy-lemma}
{\em Let $X \subset \R^n$ be a definable set, $Y \subset \R^p$ a definable $C^1$ manifold, and $f \colon X \to Y$ a proper definable $C^1$ map. 
(Furthermore, assume that $f$ satisfies the boundedness condition (II) in Remark \ref{boundedness-conditions} unless the category is $\X_0$.)
Suppose that $X$ admits a finite definable $C^1$ Whitney stratification $\{X_\alpha\}$ and $f|_{X_\alpha} \colon X_\alpha \to Y$ is a $C^1$ submersion onto $Y$ for each $\alpha$. Then, $f$ is a definable $C^0$ locally trivial fibration whose restriction $f|_{X_\alpha}$ to each $X_\alpha$ is a definable $C^1$ locally trivial fibration.}
\end{theorem}

\subsubsection{Approximation theorem}

In our later argument, we will focus on one stratum (top dimensional stratum) of $X$. Using the following Theorem \ref{approximation-thm}, we can find a definable trivialization map whose regularity is high enough as we want near the $C^1$ trivialization map obtained by the above Theorem \ref{isotopy-lemma}.

As its preparation, let us shortly recall about spaces of definable maps and their topology (the following notation refers to \cite{Val}).
Let $\cal{D}^1(X, Y)$ denote the space of definable $C^1$ maps between definable $C^1$ manifolds $X$ and $Y$.
For $f \in \cal{D}^1(X, \R)$, consider the derivative $df \colon X \to \R^{\dim X}$ and define
\[|f|_1 \coloneqq |f| + |df| \colon X \to \R,\]
where $|\hy|$ is the Euclidean norm.
Also, for a definable function $\epsilon \colon X \to (0, \infty)$ (in case of $\frak{X}$, require that $\epsilon$ is bounded), define
\[\cal{U}_\epsilon^1(f) \coloneqq \set{g \in \cal{D}^1(X, \R)}{|f - g|_1 < \epsilon}.\]
Such $\cal{U}_\epsilon^1(f)$ produce the {\em definable $C^1$ topology} on $\cal{D}^1 (X, \R)$. 
It is similar for $\cal{D}^1 (X, Y)$. 

\begin{theorem}[approximation theorem for definable maps {\cite[(II.5.2)]{Shi_G}}]\label{approximation-thm}
{\em Let $X, Y$ be definable $C^r$ manifolds and $f \colon X \to Y$ be a definable $C^1$ map. 
(Furthermore, assume that $f$ satisfies the boundedness condition (I) in Remark \ref{boundedness-conditions} unless the category is $\X_0$.)
Then, $f$ can be approximated by a definable $C^r$ map arbitrary closely in the definable $C^1$ topology.}
\end{theorem}

\begin{proposition}\label{Cr-trivial}
{\em Suppose the same setup as in Theorem \ref{isotopy-lemma}. Assume that $f: X\to Y$ and the stratification $\{X_\alpha\}$ are of class $C^r$. Choose one stratum $X_0$. 
Then, by modifying the trivialization map, $f|_{X_0} \colon X_0 \to Y$ becomes a definable $C^r$ locally trivial fibration.}
\end{proposition}

\begin{Proof}
In this proof, let $f$ mean $f|_{X_0} \colon X_\alpha \to Y$ for short.
Take a definable $C^1$ locally trivial fibration structure of $f$ using Theorem \ref{isotopy-lemma}.
Let $y \in Y$ be a point, $V$ an open neighborhood of $y$ in $Y$, and $(f, \Phi) \colon f^{-1}(V) \to V \times f^{-1}(y)$ the definable $C^1$ trivialization map of $f$ over $V$. 
Put $S = \cal{D}^1(f^{-1}(V), V \times f^{-1}(y))$ with the definable $C^1$ topology. 
Note that the subset of all diffeomorphisms is open in $S$ (see \cite[(II.5.3)]{Shi_G} for the detail). 
Hence, applying Theorem \ref{approximation-thm} to the map $\Phi \colon f^{-1}(V) \to f^{-1}(y)$, we can find a definable $C^r$ map $\phi \colon f^{-1}(V) \to f^{-1}(y)$ near $\Phi$ such that $(f, \phi)$ is a definable $C^r$ diffeomorphism. $\qed$
\end{Proof}

\begin{remark}[regularity $(r = \omega)$ in semialgebraic/subanalytic case]\label{remark2}
Theorem \ref{approximation-thm} (and hence, Proposition \ref{Cr-trivial}) hold even if we replace ``definable $C^r$'' as ``semialgebraic $C^\omega$ (that is Nash)'' and ``subanalytic $C^\omega$'' \cite{Shi_N, Whi_0}.
\end{remark}

\subsection{Definable Lie groupoids}

We recall some notions of Lie groupoid theory.
For a detailed account of basic Lie groupoid theory, see \cite{Mac}.

\begin{definition}[definable $C^r$ Lie groupoid]
A {\em definable $C^r$ Lie groupoid} $\cal{G} \toto M$ consists of the following data with conditions.
First, the data are
\begin{description}
\item[space of arrows] a definable $C^r$ manifold $\cal{G}$;
\item[space of objects] a definable $C^r$ manifold $M$;
\item[source map and target map] two surjective definable $C^r$ submersions $s, t \colon \cal{G} \to M$;
\item[composition map] a definable $C^r$ map $c \colon \cal{G}^{(2)} \to \cal{G}$, $g \cdot h \coloneqq c(g, h)$, where $\cal{G}^{(2)} = \set{(g, h) \in \cal{G} \times \cal{G}}{s(g) = t(h)}$;
\item[unit section] a definable $C^r$ embedding $u \colon M \to \cal{G}$;
\item[inverse map] a definable $C^r$ diffeomorphism $i \colon \cal{G} \to \cal{G}$, $g^{-1} \coloneqq i(g)$.
\end{description}
Second, the required conditions are that for all $g, h, k \in \cal{G}$ and $m \in M$, the following properties hold whenever they are defined:
\begin{description}
\item[composition] $(g \cdot h) \cdot k = g \cdot (h \cdot k)$, $s(g \cdot h) = s(h)$, $t(g \cdot h) = t(g)$;
\item[unit] $s(u(m)) = t(u(m)) = m$;
\item[inverse] $g^{-1} \cdot g = u(s(g))$, $g \cdot g^{-1} = u(t(g))$.
\end{description}
\end{definition}

Below we simply say ``definable'' Lie groupoids/manifolds/maps/locally trivial fibrations to mean ``definable $C^r$'', unless specifically mentioned.

Let $\cal{G} \toto M$ be a definable Lie groupoid.

\begin{definition}[saturation and orbit]
For any subset $S \subset M$, we define the {\em ($\cal{G}$-)saturation of $S$} as
\[\cal{G}.S \coloneqq \set{t(g) \in M}{g \in \cal{G} \text{ and } s(g) \in S} = t(s^{-1}(S)).\]
When $S = \{m\}$, its saturation $\cal{G}.S = \cal{G}.m$ of $S$ is called the {\em ($\cal{G}$-)orbit of the point $m$}.
A subset $S$ of $M$ is said to be {\em $\cal{G}$-invariant} if $\cal{G}.S = S$.
\end{definition}

Notice that if $S$ is a definable subset of $M$ (and in case of $\X$, if $s$ and $t$ satisfy the  boundedness condition), then $\cal{G}.S$ of $S$ is also definable.

If a definable submanifold $X$ of $M$ is $\cal{G}$-invariant, then the restriction $\cal{G}|_X \toto X$ is also a definable Lie groupoid, where $\cal{G}|_X \coloneqq s^{-1}(X) = t^{-1}(X)$. 

Definable Lie groupoids have a {\em partial} homogeneity around each orbits because the following lemma holds. 

\begin{lemma}
{\em Let $g \in \cal{G}$.
Then, there exists a definable local {\em bisection} $\sigma(U) \subset \cal{G}$ through $g$, that is, $g$ has a definable neighborhood $U$ of $s(g)$ in $M$ and a definable map $\sigma \colon U \to \cal{G}$ such that $s \circ \sigma = \op{id}_U$ and $t \circ \sigma \colon U \to M$ is a definable diffeomorphism onto the open neighborhood $t(\sigma(U))$ of $t(g)$ in $M$.}
\end{lemma}

\begin{Proof}
See \cite[Proposition 1.4.9]{Mac} for case of plain Lie groupoid. $\qed$
\end{Proof}


At the end of this section, we construct an invariant stratification of definable Lie groupoids that follows from Proposition \ref{classical-stratification}.

\begin{proposition} \label{key-lemma}
{\em  Let $V \subset M$ be a $\cal{G}$-invariant definable closed subset and $V' \subset V$ a $\cal{G}$-invariant definable $C^r$ manifold such that $\dim V = \dim V' > \dim (V - V')$.
Then, $V$ admits a finite definable $C^r$ Whitney stratification such that $V'$ is the top stratum and each stratum is $\cal{G}$-invariant.}
\end{proposition}

\begin{Proof}
Take a  definable Whitney strtatification of $V$ with the filtration by $V_i$ as in the proof of Proposition \ref{classical-stratification}.
It suffices to check that each $V_i$ is $\cal{G}$-invariant.
Clearly, both $V_0 \coloneqq V$ and $V_1 \coloneqq V - V'$ are $\cal{G}$-invariant.
Suppose that $V_0, V_1, \dots, V_i$ are all $\cal{G}$-invariant. Since both singularity and Whitney (b)-regularity condition are invariant under local diffeomorphisms, $\Sigma_r V_i$, $B(V', V_i - \Sigma_r V_i)$, and $B(V_i - \Sigma_r V_i, V_j - \Sigma_r V_j) \ (j = 0, \dots, i-1)$ are all $\cal{G}$-invariant.
Hence $V_{i+1}$ is $\cal{G}$-invariant. $\qed$
\end{Proof}

\begin{remark}[regularity $(r = \omega)$ in semialgebraic/subanalytic case]
Proposition \ref{key-lemma} holds even if we replace ``definable $C^r$'' as ``semialgebraic $C^\omega$ (that is Nash)'' and ``subanalytic $C^\omega$.''
\end{remark}

For the following section, we define the equivalence relation
\[x \sim y \ \Longleftrightarrow \ s(g) = x \text{ and } t(g) = y \text{ for some } g \in \cal{G}\]
on $M$; we consider the quotient space $M / \cal{G} \coloneqq M / \!\! \sim$ and the quotient map $q \colon M \to M / \cal{G}$.

\section{The main result}

In this section, we state and prove the main result.

Let $r$ be a fixed integer such that $1 \le r < \infty$.
Our main result is the following theorem.

\begin{theorem}\label{main2}
{\em Let $\cal{G} \rightrightarrows M$ be a definable $C^r$ Lie groupoid. 
(Furthermore, assume that both $\cal{G}$ and $M$ are bounded unless the category is $\X_0$.)
Then, there exists a filtration
\[M = M_0 \supset M_1 \supset M_2 \supset \cdots \supset M_{d+1} = \void\]
of $M$ such that for each $i = 0, 1, \dots, d\,   (= \dim M)$,}
\begin{enumerate}
\item {\em the set $M_i$ is a $\cal{G}$-invariant definable closed subset of $M$;}
\item {\em the set $M_i - M_{i+1}$ is a definable $C^r$ manifold of codimension $i$ in $M$ (unless it is empty) and $\{M_i - M_{i+1}\}_{i=0}^d$ is a definable $C^r$ Whitney stratification of $M$;} 
\item {\em the quotient space $(M_i - M_{i+1}) / \cal{G}$ admits a $C^r$ manifold structure and the quotient map $q \colon M_i - M_{i+1} \to (M_i - M_{i+1}) / \cal{G}$ is a $C^r$ locally trivial fibration.  Moreover, the quotient manifold and the quotient map are piecewise definable.}
\end{enumerate}
\end{theorem}

Here, a {\em piecewise definable $C^r$ manifold} is a $C^r$ manifold given by an atlas $\{(U_\lambda, \phi_\lambda)\}$ such that each $\phi(U_\lambda)$ is a definable set and each $\phi_\lambda \circ \phi_\mu^{-1}$ is a definable $C^r$ map; a {\em piecewise definable $C^r$ map} is a continuous map between piecewise definable $C^r$ manifolds whose each local representation is a definable $C^r$ map.

\begin{remark}[assumption in case of $\X$]
In case of $\X$, we assume the following conditions which are equivalent to each other:
\begin{itemize}
\item $\cal{G}$ and $M$ are bounded;
\item $M$ is bounded and $s$ and $t$ satisfy the boundedness condition.
\end{itemize}
Then, we have the following properties which will be used later: 
\begin{itemize}
\item The image of $M$ and its subsets by embedding $\R^n$ into $\bold{RP}^n$ become definable;
\item For every definable subset $S \subset M$, the saturation $\cal{G}.S$ becomes definable.
\end{itemize}
\end{remark}

\begin{remark}[regularity $(r = \omega)$ in semialgebraic/subanalytic case]
Theorem \ref{main2} holds even if we replace ``definable $C^r$'' as ``semialgebraic $C^\omega$ (that is Nash)'' and ``subanalytic $C^\omega$'', because of Remark \ref{remark1}, Remark \ref{remark2}, and the following proof. Therefore, especially, we obtain Theorem \ref{main} stated in Introduction (cf. Mather \cite{Mat_IDGA} and Vassiliev \cite{Vass}). According to a recent work of Vallete--Vallete \cite{Val}, we may have a chance to improve the regularity for some restricted $\X_0$-category. 
\end{remark}

For the proof, we introduce the following notion (cf.~\cite{Mat_IDGA}).

\begin{definition}[definably smooth family of submanifolds]\label{definably-smooth}
Let $U, M$ be definable $C^r$ manifolds and for each point $u \in U$ it is assigned a definable $C^r$ submanifold $V_u$ of $M$.
Then, the family $\{V_u\}_{u \in U}$ is said to be {\em definably smooth} if
the union
\[V \coloneqq \coprod_{u \in U} u \times V_u \subset U \times M\]
is a definable $C^r$ manifold and the projection $\op{pr}_1 \colon V \to U$ is a definable $C^r$ locally trivial fibration.
\end{definition}

First we replace (3) in Theorem \ref{main2} by the following $(3')$. 

\begin{lemma}\label{main3}
{\em 
Consider the same setup as in Theorem \ref{main2}. 
Then, there is a filtration of $M$ satisfying (1) and (2) of Theorem \ref{main2} and 
}
\begin{enumerate}
\item[$(3')$] {\em The family $\{\cal{G}.x\}_{x \in M_i - M_{i+1}}$ is definably smooth.}
\end{enumerate}
\end{lemma}

We now divide the proof of Theorem \ref{main2} into two steps: we will first show that Lemma \ref{main3} implies Theorem \ref{main2}, and then we prove Lemma \ref{main3}.
An essential idea of the proof is based on a sketchy proof in Mather \cite{Mat_IDGA}.

\subsection{Reduction of Theorem \ref{main2} to Lemma \ref{main3}}

It suffices to prove that Lemma \ref{main3} implies the condition $(3)$ in Theorem \ref{main2}.

Let $\cal{G} \rightrightarrows M$ be a definable Lie groupoid. Suppose that we have a filtration by $M_i$'s as in Lemma \ref{main3} and $i$ a fixed number. Put 
\[M^i = M_i - M_{i+1}, \quad \cal{R}^i \coloneqq \coprod_{x \in M^i} x \times \cal{G}.x = \set{(x, y) \in M^i \times M^i}{x \sim y}.\]
The condition $(3')$ in Lemma \ref{main3} means that $\cal{R}^i$ is a definable manifold and the projection $\op{pr}_1 \colon \cal{R}^i \to M^i$ is a definable $C^r$ locally trivial fibration.

We now show the following four claims which imply $(3)$ in Theorem \ref{main2}. The central idea is to find a definable version of `slice theorem' (without the assumption on properness of the action). 

\begin{claim}\label{orbit-is-submfd}
{\em For every point $x \in M^i$, the orbit $\cal{G}.x$ is a closed definable submanifold of $M^i$.}
\end{claim}

{\em Proof of Claim \ref{orbit-is-submfd}.}
Take a regular point of the definable set $\cal{G}.x$. By mapping a neighborhood of the point to around other points by definable bisections through arrows in $s^{-1}(x)$, then $\cal{G}.x$ has a chart as a regular submanifold of $M^i$ around each point.

Suppose that there exists a point $y \in \op{Cl}_{M^i}(\cal{G}.x) - \cal{G}.x$. Since $\cal{G}.y \subset \op{Cl}_{M^i}(\cal{G}.x) - \cal{G}.x$ and orbits are definable, we obtain that
\[\dim \cal{G}.y \le \dim (\op{Cl}_{M^i}(\cal{G}.x) - \cal{G}.x) < \dim \cal{G}.x.\]
However, $x$ and $y$ are points of $M^i$ and $\cal{R}^i \to M^i$ is locally trivial, thus $\dim \cal{G}.y = \dim \cal{G}.x$. This makes the contradiction. 
Thus, $\cal{G}.x$ is closed in $M^i$. $\square$\\

Take an arbitrary point $a \in M^i$. Since $\cal{G}.a$ is a regular submanifold of $M^i$, there is a definable submanifold $S$ of $M^i$ passing through $a$,  such that $\cal{R}^i$ is trivialized over $S$ and $S$ intersects $\cal{G}.a$ only at the point $a$ transversely.
Hereafter, we call $S$ a {\em slice} in $M^i$ to $\cal{G}.a$ at the point $a$.
Notice that $t(\sigma(S)) \cap \cal{G}.a = \{t(\sigma(a))\}$ for any local bisection $\sigma \colon U \to \cal{G}$ on an open neighborhood $U$ of $a$ in $M^i$.

\begin{claim}\label{make-S-slice}
{\em By shrinking $S$ suitably if necessary, $\cal{G}.S$ becomes to be open in $M^i$ and $S \cap \cal{G}.x$ is empty or consists of one point for each $x \in M^i$. 
In particular, the quotient space $M^i / \cal{G}$ is Hausdorff.}
\end{claim}

{\em Proof of Claim \ref{make-S-slice}.}
Consider the restriction of the source map to $s^{-1}(S)$, 
\[\tilde{s} = s|_{s^{-1}(S)} \colon s^{-1}(S) \to M^i.\]
The intersection of $u(S)$ and $s^{-1}(a)$ at $u(a)$ is transverse in $s^{-1}(S)$, that is, $T_{u(a)} s^{-1}(S) = T_{u(a)} u(S) \oplus T_{u(a)} s^{-1}(a)$. In addition, $d\tilde{s}_{u(a)}$ maps
\[T_{u(a)} u(S) \oplus 0 \to T_a S, \quad 0 \oplus T_{u(a)} s^{-1}(a) \to T_a \cal{G}.a\]
surjectively.
Thus, $\tilde{s}$ is submersive at $u(a)$. From the implicit function theorem, $\tilde{s}$ is an open map and hence $\cal{G}.S$ becomes to be open in  $M^i$ by retaking $S$ small enough. 
Next, consider the second factor projection 
\[\rho = \op{pr}_2|_{\op{pr}_1^{-1}(S)} \colon  \op{pr}_1^{-1}(S) \to \cal{G}.S.\]
Here, by using $(3')$ of Lemma \ref{main3}, we may assume that $\op{pr}_1^{-1}(S)$ is diffeomorphic to $S \times \cal{G}.a$ such that the first factor projection to $S$ commutes with $\op{pr}_1$. 
We also write $\rho \colon S \times \cal{G}.a \to \cal{G}.S$ for short. Since
\[d\rho_{(a, a)} \colon T_a S \oplus T_a \cal{G}.a \to T_a \cal{G}.S = T_a M^i\]
is a linear isomorphism from the transversality condition, $\rho$ is a diffeomorphism on a neighborhood $U$ of $(a,a)$ in $\op{pr}_1^{-1}(S)$.
Now suppose that for any neighborhood of $a$ in $S$, there exists two points $x, y \in S$ such that $x \ne y$ and $\cal{G}.x = \cal{G}.y$. Then we can take two sequences $(x_n), (y_n)$ on $S$ such that $x_n, y_n \to a \ (n \to \infty)$ and $x_n \ne y_n$ for each $n$.
Since $(x_n, x_n), (y_n, x_n) \to (a, a) \ (n \to \infty)$, it follows that $(x_N, x_N), (y_N, x_N) \in U$ for some number $N$.
However, these two points of $U$ are mapped to $x_N \in \rho(U)$, that makes the contradiction to that $\rho$ is diffeomorphic to $U$. 
Finally, it is clear that $M^i / \cal{G}$ is Hausdorff. $\square$\\

Hereafter, we take all slices small enough as in Claim 2. 
Remark that the above $\rho$ is a definable diffeomorphism,  for it is bijective and locally diffeomorphic. 
Now we introduce a definable $C^r$ manifold structure of $M^i / \cal{G}$ by using slices.

\begin{claim}\label{orbit-space-is-mfd}
{\em The quotient space $M^i / \cal{G}$ admits a piecewise definable manifold structure.}
\end{claim}

{\em Proof of Claim \ref{orbit-space-is-mfd}.}
Take $[a] = \cal{G}.a \in M^i / \cal{G}$ and a slice $S$ at $a \in M^i$.
The restriction of the quotient map $q \colon M^i \to M^i / \cal{G}$ to $S$ is bijective since $S \cap \cal{G}.x$ has at most one point for each $x \in M^i$.
In addition, $q$ is continuous and open, and hence $q|_S$ is a homeomorphism onto its image.
Then we introduce the following chart around $[a] \in M^i / \cal{G}$ using any slice $S$ at $a$:
\[\phi \colon q(S) \to S \to  W \subset \R^s,\]
where the first arrow is $(q|_S)^{-1}$ and the second arrow is a local chart of $S$ onto an open set $W \subset \R^s$ as a definable submanifold ($s = \dim S$).
We check that $\{(q(S), \phi)\}$ forms an atlas of a piecewise definable $C^r$ manifold.
Let $S$ and $S'$ be slices with $q(S) \cap q(S') \ne \void$. For our convenience, we retake them as $q(S) = q(S')$  (i.e. $\cal{G}.S = \cal{G}.S'$).  
It suffices to show that $g \coloneqq (q|_{S'})^{-1} \circ q|_S \colon S \to S'$ is a definable $C^r$ diffeomorphism. In fact, the map $g$ is the restriction of the definable $C^r$ map 
\[\tilde{g} \coloneqq \op{pr}_1 \circ \rho^{-1} \colon \cal{G}.S' \xto{\isom} S' \times \cal{G}.a \to S'\]
to $S$, where $\Phi$ is a trivialization map of $\op{pr}_1 \colon \cal{R}^i \to M^i$ on $S'$.
Therefore, $g = \tilde{g}|_S$ is definable, of class $C^r$, and non-singular. $\square$\\

\begin{claim}\label{quot-map-is-loc-triv}
{\em The quotient map $q \colon M^i \to M^i / \cal{G}$ is a piecewise definable locally trivial fibration.}
\end{claim}

{\em Proof of Claim \ref{quot-map-is-loc-triv}.}
For any chart $(S, \phi)$ of the above, $q|_{\cal{G}.S} \colon \cal{G}.S = q^{-1}(q(S)) \to q(S)$ is locally expressed as $\op{pr}_1 \colon S \times \cal{G}.a \to S$ via $\phi$. $\square$\\

This completes the proof of $(3)$ in Theorem \ref{main2}. $\qed$

\subsection{Proof of Lemma \ref{main3}}

Let $\cal{G} \toto M \, (\subset \R^n)$ be a definable Lie groupoid.
We prove Lemma \ref{main3} by the induction on the codimension $i$. 
Assume that we have $M = M_0 \supset M_1 \supset \dots \supset M_i$ such that for every $j=1, \cdots, i$, it holds that 
\begin{enumerate}
\item[$(1)_{j}$] the set $M_{j}$ is a $\cal{G}$-invariant definable closed subset of $M$ and $\op{codim} \ M_j \ge j$;
\item[$(2)_{j}$] the set $M_{j-1} - M_{j}$ is a definable $C^r$ manifold of codimension $j-1$ in $M$ and $\{M_k - M_{k+1}\}_{k=0}^{j-1}$ is a definable $C^r$ Whitney stratification of $M - M_{j}$;
\item[$(3)_{j}$] the family $\{\cal{G}.x\}_{x \in M_{j-1} - M_j}$ is definably smooth.
\end{enumerate}
It suffices to find a subset $M_{i+1} \subset M_i$ such that $(1)_{i+1}$, $(2)_{i+1}$, and $(3)_{i+1}$ hold. \\

If $\op{codim} \ M_i > i$, then it suffices to take $M_{i+1} \coloneqq M_i$. So we assume that $\op{codim} \ M_i = i$.
We first put
\[\cal{R}_i := \coprod_{x \in M_i} x \times \cal{G}.x = \set{(x, y) \in M_i \times M_i}{x \sim y}.\]
Note that the first factor projection $\op{pr}_1 \colon \cal{R}_i \to M_i$ may not be proper.
To compactify fibers of $\cal{R}_i$, we embed $\R^n$ into $\bold{RP}^n$ by $(x_1, \dots, x_n) \mapsto [1 : x_1 : \dots : x_n]$, where $M \subset \R^n$. Further, embed $\bold{RP}^n$ into some $\R^N$, for example, by the following semialgebraic map:
\[\bold{RP}^n \to \R^{n+1+\frac{n(n+1)}{2}}, \quad [x_0 \colon x_1 \colon \dots \colon x_n] \mapsto \left( \frac{x_i x_j}{\sum_{k=0}^n x_k^2} \right)_{0 \le i \le j \le n}.\]
We set 
\[\overline{M_i} = \op{Cl}_{\bold{RP}^n} M_i \subset \bold{RP}^n, \quad \overline{\cal{R}_i} = \op{Cl}_{M_i \times \overline{M_i}} \cal{R}_i \subset M_i \times \overline{M_i}.\]
Then $\overline{M_i}$ is compact and both $\overline{M_i}$ and $\overline{\cal{R}_i}$ are definable.

We now consider a definable Lie groupoid $\cal{G} \toto M \times \R^N$ with the action on the first factor; the sets $\cal{R}_i$ and $\overline{\cal{R}_i}$ are $\cal{G}$-invariant.
Hence, by using Proposition \ref{key-lemma}, we have a $\cal{G}$-invariant definable Whitney stratification $\cal{S} = \{S_\alpha\}$ of $\overline{\cal{R}_i}$ such that the top stratum is $\cal{R}_i - \Sigma_r \cal{R}_i$. 
We set
\[\displaystyle M_{i+1} := \op{Cl}_{M_i} \left( \Sigma_r M_i \cup \bigcup_{j = 0}^{i-1} B(M_j - M_{j+1}, M_i - \Sigma_r M_i) \cup \bigcup_\alpha \pi(C(\pi|_{S_\alpha})) \cup \pi(\Sigma_r \cal{R}_i) \right),\]
where $C(h)$ denotes the critical point set of a $C^r$ map $h$. 

We show that $M_{i+1}$ satisfies $(1)_{i+1}$, $(2)_{i+1}$, and $(3)_{i+1}$.\\

\noindent 
$(1)_{i+1}: \quad$ It is obvious that $\Sigma_r M_i$, $B(M_j - M_{j+1}, M_i - \Sigma_r M_i)$, $C(\pi|_{S_\alpha})$, and $\Sigma_r \cal{R}_i$ are $\cal{G}$-invariant and definable.
Hence, we see that $\pi(C(\pi|_{S_\alpha}))$ and $\pi(\Sigma_r \cal{R}_i)$ are $\cal{G}$-invariant and definable (remember that $\pi$ is $\cal{G}$-equivariant).
Thus, $M_{i+1}$ is $\cal{G}$-invariant and definable.
Next, $\Sigma_r M_i$ and $B(M_j - M_{j+1}, M_i - \Sigma_r M_i)$ are nowhere dense in $M_i$ for dimensional reason, and $\pi(C(\pi|_{S_\alpha}))$ is nowhere dense in $M_i$ from Sard's theorem.
Moreover, we see that $\pi(\Sigma_r \cal{R}_i)$ is nowhere dense as follows. 
Suppose that there exists a non-empty open subset $U$ of $M_i$ included in $\pi(\Sigma_r \cal{R}_i)$. Then, $\pi^{-1}(U)$ is also an open set in $\cal{R}_i$ such that
\[\void \ne \pi^{-1}(U) \subset \pi^{-1}(\pi(\Sigma_r \cal{R}_i)) = \Sigma_r \cal{R}_i.\]
This makes the contradiction to that $\Sigma_r \cal{R}_i$ is nowhere dense in $\cal{R}_i$.
Consequently, $M_{i+1}$ is also nowhere dense in $M_i$, and hence, we have that $\op{codim} \ M_{i+1} > \op{codim} \ M_i$.
\\

\noindent 
$(2)_{i+1}: \quad$ Note that $M_i - M_{i+1}$ is open in $M_i - \Sigma_r M_i$, thus $M_i - M_{i+1}$ is a definable manifold of codimension $i$.
Moreover, for each $j = 0, 1, \dots, i-1$, since $M_{i+1}$ contains the bad set $B(M_j - M_{j+1}, M_i - \Sigma_r M_i)$, the pair $(M_j - M_{j+1}, M_i - M_{i+1})$ satisfies the Whitney (b)-regularity condition.
Thus, $\{M_j - M_{j+1}\}_{j = 0}^i$ is a definable Whitney stratification of $M - M_{i+1}$.
\\

\noindent 
$(3)_{i+1}: \quad$ Let $X$ denote the preimage $\pi^{-1}(M_i - M_{i+1}) \subset \overline{\cal{R}_i}$.
Then, $X$ has the induced stratification from $\cal{S}=\{S_\alpha\}$ of $\overline{\cal{R}_i}$ and we will apply Theorem \ref{Cr-trivial} to the definable map $\pi = \pi|_X \colon X \to M_i - M_{i+1}$.
To do this, we check that $\pi$ satisfies the assumption in Theorem \ref{Cr-trivial}.
It is obvious that $\pi$ is proper (for we compactified fibers of $\cal{R}_i$). 
Moreover, for each stratum $S_\alpha \in \cal{S}$, we see that $\pi$ is submersive on $S_\alpha \cap X$ (since $\pi(C(\pi|_{S_\alpha}))$ is included in $M_{i+1}$). 
Here, notice that the top dimensional stratum of $X$ is 
\[X_0 \coloneqq (\cal{R}_i - \Sigma_r \cal{R}_i) \cap X = \cal{R}_i \cap \pi^{-1}(M_i - M_{i+1}),\]
for $\Sigma_r\cal{R}_i \cap \pi^{-1}(M_i - M_{i+1})=\emptyset$ by the definition of $M_{i+1}$. 
Now applying Proposition \ref{Cr-trivial} to the map $\pi \colon  X \to M_i - M_{i+1}$ and the stratum $X_0$, we see that
\[\pi|_{X_0} \colon  X_0 \to M_i - M_{i+1}\]
is a definable $C^r$ locally trivial fibration.
That is equivalent to that the family $\{\cal{G}.y\}_{y \in M_i - M_{i+1}}$ is definably smooth.  This completes the proof. $\qed$

\section*{Acknowledgement}
The author wishes to thank Dr.~Toru Ohmoto, his supervisor, for guiding him to this subject and for instructions and discussions, and Dr.~Satoshi Koike for valuable comments on our results.

\end{document}